\documentclass[11pt,reqno,a4paper]{amsart}
\usepackage[all]{xy}
\usepackage[margin=28mm]{geometry}
\usepackage[T1]{fontenc}
\usepackage{lmodern}
\usepackage{amsmath,amssymb,amsthm}
\usepackage[hidelinks]{hyperref}
\numberwithin{equation}{section}
\newtheorem{theorem}{Theorem}[section]
\newtheorem{proposition}[theorem]{Proposition}
\newtheorem{lemma}[theorem]{Lemma}

\theoremstyle{definition}
\newtheorem{question}[theorem]{Question}
\newcommand{\pd}{\operatorname{pd}}
\newcommand{\id}{\operatorname{id}}
\newcommand{\Ext}{\operatorname{Ext}}
\newcommand{\Tor}{\operatorname{Tor}}
\newcommand{\Hom}{\operatorname{Hom}}
\newcommand{\Mod}{\operatorname{Mod}}
\newcommand{\End}{\operatorname{End}}
\newcommand{\Soc}{\operatorname{Soc}}
\newcommand{\op}{\mathrm{op}}
\newcommand{\Ker}{\operatorname{Ker}}
\newcommand{\im}{\operatorname{Im}}
\title[The socle of the last term in a minimal injective resolution II]
{The socle of the last term in a minimal injective resolution II}
\author[Z. Huang]{Zhaoyong Huang}
\author[T. Zhao]{Tiwei Zhao$^*$}\thanks{$^*$Corresponding author.}
\address[Z. Huang]{Department of Mathematics, Nanjing University,
Nanjing 210093, Jiangsu Province, P. R. China}
\email{huangzy@nju.edu.cn}
\address[T. Zhao]{School of Artificial Intelligence, Jianghan University,
Wuhan 430056, Hubei Province, P. R. China}
\email{tiweizhao@jhun.edu.cn}
\date{}
\hypersetup{pdftitle={The socle of the last term in a minimal injective resolution II}}
\begin{document}
\begin{abstract}
Let $R$ be a Gorenstein ring with left and right self-injective dimensions $n$.
We prove that the last term of a minimal injective resolution of $_RR$ has an
essential socle. It answers an open question of Iwanaga and Sato. 
As a consequence, we get that for any left $R$-module $M$ with injective dimension $n$,
the last term of a minimal injective resolution of $M$ also has an essential socle.
\end{abstract}
\maketitle

\section{Introduction}\label{sec:introduction}

Let $R$ be an arbitrary associative ring. For a left $R$-module $M$, we use 
$$0\to M\to E^0(M)\to E^1(M)\to\cdots\to E^i(M)\to\cdots$$
to denote a minimal injective resolution of $M$, and write
$\Omega^{-i}(M):=\Ker(E^i(M)\to E^{i+1}(M))$ for any $i\geq 0$.
We identify $M$ with its image in $E^0(M)$, so that $\Omega^0(M)=M$. 
In this minimal resolution, each $E^i(M)$ is the injective envelope of $\Omega^{-i}(M)$;
in particular, $E^0(M)$ is the injective envelope of $M$.
Recall that a left and right Noetherian ring $R$ is {\it Gorenstein} if its left and right 
self-injective dimensions are finite. 

Let $R$ be a commutative Noetherian ring. Then there is a one-to-one correspondence 
between the prime ideals of $R$ and the isomorphism classes of indecomposable injective $R$-modules given by
$P\leftrightarrow E^0(R/P)$ with $P$ a prime ideal of $R$ \cite[Proposition 3.1]{M}. 
Let $R$ be a commutative Gorenstein ring with self-injective dimension $n$. Then
the Krull dimension of $R$ is $n$, and the fundamental theorem in \cite{B} states that 
an indecomposable injective $R$-module $E^0(R/P)$ 
with $P$ a prime ideal of $R$ embeds in $E^i(R)$ 
if and only if the height of $P$ is $i$. Thus the last term $E^n(R)$ has 
an essential socle. It is reasonable to ask whether this result holds true
in the non-commutative case.

\begin{question} (\cite[p.62, Question]{IS})\label{ques-1.2}
Let $R$ be a Gorenstein ring with left and right self-injective dimensions $n$. 
Is the socle of $E^n(R)$ non-zero? Moreover, is this socle essential in $E^n(R)$?
\end{question}

It is trivial that an affirmative answer to the second question entails that to the first.
For the first question, the answer is affirmative when one of the following conditions is satisfied:
\begin{enumerate}
\item[{\rm (1)}] the left and right self-injective dimensions of $R$ are at most two (\cite[Theorem~4.5]{Ho}); 
\item[{\rm (2)}] $R$ is an Auslander--Gorenstein ring (\cite[Proposition~1.1]{FI});
\item[{\rm (3)}] either the global dimension of $R$ is $n$ or $R$ is fully bounded (\cite[Proposition~8]{IS}); 
\item[{\rm (4)}] $R$ is an Auslander $(n-2)$-Gorenstein ring (\cite[Corollary~3.2]{HW}).
\end{enumerate}
For the second question, the answer is affirmative if $R$ is an Auslander--Gorenstein ring (\cite[Theorem~6]{IS}).
The module version of this result is as follows: 
for an Auslander--Gorenstein ring $R$ with left and right self-injective dimensions $n$,
if $M$ is a left $R$-module with injective dimension $n$, then $E^n(M)$ has an 
essential socle (\cite[Theorem~2]{IM}). In addition, for an algebra $R$ that is module-finite over a commutative Noetherian ring,
if $M$ is a non-zero finitely generated left $R$-module of injective dimension $n$,
then $E^n(M)$ has an essential socle (\cite[Corollary~3.5(2)]{GN}).
Let $R$ be a left Noetherian ring and let $M$ be a Wakamatsu tilting left $R$-module with $S:=\End(_RM)^{op}$
a right Noetherian ring. If $M$ is a Gorenstein module with left $R$-right $S$-injective dimensions $n$, then
$E^n(M)$ has an essential socle (\cite[Definition~3.6 and Corollary 3.10]{SZH}).

In this paper, we prove the following result. 

\begin{theorem}\label{thm-1.2}
{\rm (Theorems \ref{thm-3.2} and \ref{thm-3.3})}
Let $R$ be a Gorenstein ring with left and right self-injective dimensions $n$. 
Then it holds that
\begin{enumerate}
\item[{\rm (1)}] 
$E^n(R)$  has an essential socle.
\item[{\rm (2)}] 
$E^n(M)$ has an essential socle for any non-zero left $R$-module $M$ with injective dimension $n$. 
\end{enumerate}
\end{theorem}

In this theorem, assertion (1) answers the second question 
(and hence the first one) in Question \ref{ques-1.2} affirmatively, and assertion (2) is the module version of (1).

\section{Preliminaries}\label{sec:preliminaries}

Throughout this paper, $R$ is an associative ring with identity,
all modules are unital and $\Mod R$ is the category of left $R$-modules.
Let $M$ be a module in $\Mod R$. We use $\pd_RM$ and $\id_RM$ to denote the projective
and injective dimensions of $M$, respectively. The {\it socle} $\Soc(M)$ is 
the sum of the simple submodules of $M$. Recall that a submodule $L$ of $M$ is 
{\it essential} if $M'\cap L\neq 0$ for any non-zero submodule $M'$
of $M$. Thus $\Soc(M)$ is essential in $M$ if and only if any non-zero submodule of $M$ 
contains a simple submodule.

\begin{lemma}\label{lem-2.1}
Let
$$M_1\buildrel{f_1}\over\longrightarrow M_2\buildrel{f_2}\over\longrightarrow M_3$$
be a sequence in $\Mod R$. If $f_1$ is epic,
then $\Ker f_2=f_1(\Ker f_2f_1)$. 
\end{lemma}

\begin{proof}
It is trivial that $f_1(\Ker f_2f_1)\subseteq\Ker f_2$. Now let $x_2\in\Ker f_2$.
Since $f_1$ is epic, there exists $x_1\in M_1$ such that $x_2=f_1(x_1)$.
So $f_2f_1(x_1)=f_2(x_2)=0$, and hence $x_1\in\Ker f_2f_1$. It follows that
$x_2=f_1(x_1)\in f_1(\Ker f_2f_1)$, and thus $\Ker f_2\subseteq f_1(\Ker f_2f_1)$. 
The proof is finished.
\end{proof}

For a submodule $X$ of a module $Y\in\Mod R$, we use $\lambda^X_Y:X\hookrightarrow Y$
to denote the inclusion.

\begin{lemma}\label{lem-2.2}
If $\id_RR=n$, then for any submodule $X$ of a module $Y\in\Mod R$, the homomorphism
$\Ext_R^n(\lambda^X_Y,R)$ in $\Mod R^{op}$ is epic.
\end{lemma}

\begin{proof}
The exact sequence
$$0\to X\buildrel{\lambda^X_Y}\over\longrightarrow Y\to  Y/X\to 0$$
in $\Mod R$ induces an exact sequence
$$\Ext_R^n(Y,R)\buildrel{\Ext_R^n(\lambda^X_Y,R)}\over\longrightarrow\Ext_R^n(X,R)\to\Ext_R^{n+1}(Y/X,R)$$
in $\Mod R^{op}$. Since $\id_RR=n$, we have $\Ext_R^{n+1}(Y/X,R)=0$, and the assertion follows.
\end{proof}

We need the following result.

\begin{lemma}\label{lem-2.3}
Let $0\neq M\in\Mod R$ and $n\geq 0$. It holds that
\begin{enumerate}
\item[{\rm (1)}] $\Ext_R^n(K,M)\neq 0$ for any non-zero submodule $K$ of $\Omega^{-n}(M)$;
\item[{\rm (2)}] if $\id_RM=n$, then $\Ext_R^n(K,M)\neq 0$ for any non-zero submodule $K$ of $E^{n}(M)$.
\end{enumerate}
\end{lemma}

\begin{proof}
(1) Let $K$ be a non-zero submodule of $\Omega^{-n}(M)$.
If $n=0$, then the inclusion $\lambda_M^K:K\hookrightarrow M$ is a non-zero element of
$\Hom_R(K,M)=\Ext_R^0(K,M)$ since $\Omega^0(M)=M$.
Thus the assertion holds in this case.
Now suppose $n\geq 1$. Consider the following diagram with exact row:
$$\xymatrix@=0.5cm{
& & & K\ar@{^(->}[d]^{\lambda_{\Omega^{-n}(M)}^K} &\\
0\ar[r] & \Omega^{-(n-1)}(M)\ar[r] & E^{n-1}(M)\ar[r]^{\ \pi^n} & \Omega^{-n}(M)\ar[r] & 0.}$$
We claim that there is not a homomorphism $f:K\to E^{n-1}(M)$ such that $\lambda_{\Omega^{-n}(M)}^K=\pi^nf$.
Otherwise, if there exists a homomorphism $f:K\to E^{n-1}(M)$ such that $\lambda_{\Omega^{-n}(M)}^K=\pi^nf$,
then $f$ is monic, and hence $\im f\cong K\neq 0$.
Thus $\im f\cap\Omega^{-(n-1)}(M)\neq 0$ as $\Omega^{-(n-1)}(M)$ is essential in $E^{n-1}(M)$.
So there exists $x\in K$ such that $0\neq f(x)\in\Omega^{-(n-1)}(M)$, and hence
$$0\neq x=\lambda_{\Omega^{-n}(M)}^K(x)=\pi^nf(x)=0,$$
which is a contradiction. The claim is proved. It follows that $\Hom_R(K,\pi^n)$ is not epic. Then the 
following exact sequence
$$0\to\Hom_R(K,\Omega^{-(n-1)}(M))\to\Hom_R(K,E^{n-1}(M))$$
$$\buildrel{\Hom_R(K,\pi^n)}\over\longrightarrow\Hom_R(K,\Omega^{-n}(M))
\to\Ext^1_R(K,\Omega^{-(n-1)}(M))\to0$$
yields $\Ext^1_R(K,\Omega^{-(n-1)}(M))\neq 0$. Thus $\Ext_R^n(K,M)\cong\Ext^1_R(K,\Omega^{-(n-1)}(M))\neq 0$
by dimension shifting.

(2) If $\id_RM=n$, then $\Omega^{-n}(M)=E^{n}(M)$. Now the assertion follows from (1).
\end{proof}

\section{Main results}
The following result is useful in the sequel.

\begin{proposition}\label{prop-3.1}
Let $N\in\Mod R^{op}$ with $\pd_{R^{\op}}N=n\geq 0$. Then for any $M\in\Mod R$
and any family $\{M_i\mid i\in I\}$ of submodules of $M$, 
the inclusion-induced maps $\Tor_n^R(N,M_i)\to\Tor_n^R(N,M)$ are monic, and
$$\Tor_n^R\left(N,\cap_{i\in I}M_i\right)=\cap_{i\in I}\Tor_n^R(N,M_i)\eqno{(3.1)}$$
as subgroups of $\Tor_n^R(N,M)$.
\end{proposition}

\begin{proof}
Let $P\in\Mod R^{op}$ be projective. Then for any submodule $X$ of $M$, the homomorphism
$$1_P\otimes\lambda_M^X:P\otimes_R X\longrightarrow P\otimes_R M$$
of Abelian groups is monic. We identify $P\otimes_R X$ with its image in $P\otimes_RM$.

We claim that
$$P\otimes_R(\cap_{i\in I}M_i)=\cap_{i\in I}(P\otimes_R M_i)$$
as subgroups of $P\otimes_R M$.
For any free right module $R^{(J)}$ with $J$ a set, we have $R^{(J)}\otimes_R X\cong X^{(J)}$. In addition, 
we have $(\cap_{i\in I}M_i)^{(J)}=\cap_{i\in I}M_i^{(J)}$ by a direct verification. It follows that
$$R^{(J)}\otimes_R(\cap_{i\in I}M_i)\cong(\cap_{i\in I}M_i)^{(J)}
=\cap_{i\in I}M_i^{(J)}\cong\cap_{i\in I}(R^{(J)}\otimes_R M_i).\eqno{(3.2)}$$
Now choose a projective right module $P'$ such that
$P\oplus P'=R^{(J)}$ with $J$ a set. For any submodule $X$ of $M$, we have
$$R^{(J)}\otimes_R X=(P\otimes_R X)\oplus(P'\otimes_R X)\quad \text{and}\quad
(P\otimes_RM)\cap(R^{(J)}\otimes_R X)=P\otimes_R X\eqno{(3.3)}$$
as subgroups of $R^{(J)}\otimes_RM$. It follows that
\begin{align*}
 P\otimes_R(\cap_{i\in I}M_i)
 &=(P\otimes_RM)\cap(R^{(J)}\otimes_R(\cap_{i\in I}M_i))\quad\text{(by (3.3))}\\
 &=(P\otimes_R M)\cap(\cap_{i\in I}(R^{(J)}\otimes_RM_i))\quad\text{(by (3.2))}\\
 &=\cap_{i\in I}((P\otimes_RM)\cap(R^{(J)}\otimes_RM_i))\\
 &=\cap_{i\in I}(P\otimes_RM_i).\quad\text{(by (3.3))}
\end{align*}
The claim is proved.

Let $X$ be any submodule of $M$. Since $\pd_{R^{\op}}N=n$, it is easy to see that 
$\Tor_n^R(N,\lambda_M^X):\Tor_n^R(N,X)\to\Tor_n^R(N,M)$ is monic. 
Thus the inclusion-induced map $\Tor_n^R(N,\lambda_M^{M_i}):\Tor_n^R(N,M_i)\to\Tor_n^R(N,M)$ 
is monic for any $i\in I$.

In the following, we prove that the equality (3.1) holds true. The case for $n=0$ follows from
the above claim. Now suppose that $n\geq 1$ and 
$$0\to P_n\xrightarrow{d_n}P_{n-1}\to\cdots\to P_0\to N\to 0$$
is a projective resolution of $N$ in $\Mod R^{op}$.
We identify $\Tor_n^R(N,X)$ with $\im\Tor_n^R(N,\lambda_M^X)$.
Consider the following commutative diagram with exact rows:
$$\xymatrix@C=1.3cm@R=0.9cm{
 0\ar[r] & \Tor_n^R(N,X)\ar[r]\ar@{^{(}->}[d]_{\Tor_n^R(N,\lambda_M^X)}
 & P_n\otimes_R X\ar[r]^{d_n\otimes 1_X}
   \ar@{^{(}->}[d]^{1_{P_n}\otimes\lambda_M^X}
 & P_{n-1}\otimes_R X
   \ar@{^{(}->}[d]^{1_{P_{n-1}}\otimes\lambda_M^X}\\
 0\ar[r] & \Tor_n^R(N,M)\ar[r]
 & P_n\otimes_R M\ar[r]^{d_n\otimes 1_M}
 & P_{n-1}\otimes_RM.}$$
For any $z\in P_n\otimes_R X$, we have
\begin{align*}
& z=(1_{P_n}\otimes\lambda_M^X)(z)\in\Tor_n^R(N,M)\\
\Longleftrightarrow &
(1_{P_{n-1}}\otimes\lambda_M^X)(d_n\otimes 1_X)(z)=(d_n\otimes 1_M)(1_{P_n}\otimes\lambda_M^X)(z)=0\\
\Longleftrightarrow & (d_n\otimes 1_X)(z)=0\\
\Longleftrightarrow & z\in\Tor_n^R(N,X).
\end{align*}
Consequently, we conclude that 
$$\Tor_n^R(N,X)=\Tor_n^R(N,M)\cap(P_n\otimes_R X)\eqno{(3.4)}$$
as subgroups of $P_n\otimes_RM$. It follows that
\begin{align*}
&\Tor_n^R(N,\cap_{i\in I}M_i)\\
=&\Tor_n^R(N,M)\cap(P_n\otimes_R(\cap_{i\in I}M_i))\quad\text{(by (3.4))}\\
=&\Tor_n^R(N,M)\cap(\cap_{i\in I}(P_n\otimes_RM_i))\quad\text{(by the above claim)}\\
=&\cap_{i\in I}(\Tor_n^R(N,M)\cap(P_n\otimes_RM_i))\\
=&\cap_{i\in I}\Tor_n^R(N,M_i)\quad\text{(by (3.4))}
\end{align*}
as subgroups of $\Tor_n^R(N,M)$.
\end{proof}

We are now in a position to prove the main result.

\begin{theorem}\label{thm-3.2}
Let $R$ be a Gorenstein ring with $\id_RR=\id_{R^{op}}R=n$. 
Then $E^n(R)$ has an essential socle. 
\end{theorem}

\begin{proof}

Let $W$ be a non-zero finitely generated submodule of $E^n(R)$.
Consider the set
$$\mathcal{C}:=\{\Ker\Ext_R^n(\lambda_W^V,R)\mid V\ \text{is a non-zero submodule of}\ W\}.$$
It is a non-empty family of submodules of $\Ext_R^n(W,R)$. Notice that $\Ext_R^n(W,R)$ is a Noetherian right $R$-module by \cite[Lemma 2]{HC}, 
so $\mathcal{C}$ has a maximal element, denoted by $\Ker\Ext_R^n(\lambda_W^{V_0},R)$.

Let $U$ be a non-zero submodule of $V_0$. Then $\lambda_W^U=\lambda_W^{V_0}\lambda^U_{V_0}$, and so
$$\Ext_R^n(\lambda_W^U,R)=\Ext_R^n(\lambda^U_{V_0},R)\Ext_R^n(\lambda_W^{V_0},R).\eqno{(3.5)}$$
It implies $\Ker\Ext_R^n(\lambda_W^{V_0},R)\subseteq\Ker\Ext_R^n(\lambda_W^U,R)$. By the maximality of $\Ker\Ext_R^n(\lambda_W^{V_0},R)$,
we have $$\Ker\Ext_R^n(\lambda_W^{V_0},R)=\Ker\Ext_R^n(\lambda_W^U,R).\eqno{(3.6)}$$
By Lemma \ref{lem-2.2}, we have that $\Ext_R^n(\lambda_W^{V_0},R)$ is epic. It follows that
\begin{align*}
&\Ker\Ext_R^n(\lambda_{V_0}^U,R)\\
=&\Ext_R^n(\lambda_W^{V_0},R)(\Ker(\Ext_R^n(\lambda^U_{V_0},R)\Ext_R^n(\lambda_W^{V_0},R)))\quad\text{(by Lemma \ref{lem-2.1})}\\
=&\Ext_R^n(\lambda_W^{V_0},R)(\Ker\Ext_R^n(\lambda_W^U,R))\quad\text{(by (3.5))}\\
=&\Ext_R^n(\lambda_W^{V_0},R)(\Ker\Ext_R^n(\lambda_W^{V_0},R))=0.\quad\text{(by (3.6))}
\end{align*}
So $\Ext_R^n(\lambda_{V_0}^U,R)$ is monic, and hence it is an isomorphism by Lemma \ref{lem-2.2}. 

Let $Q$ be an injective cogenerator for $\Mod R^{op}$. 
By \cite[Theorem 3.2.13]{EJ}, we have the following commutative diagram:
$$\xymatrix@C=3cm@R=0.9cm{
 \Tor_n^R(Q,U)\ar[r]^{\Tor_n^R(Q,\lambda_{V_0}^U)}\ar[d]_{\cong}
 & \Tor_n^R(Q,V_0)\ar[d]^{\cong}\\
\Hom_{R^{\op}}(\Ext_R^n(U,R),Q)
 \ar[r]^{\Hom_{R^{\op}}(\Ext_R^n(\lambda_{V_0}^U,R),Q)\ }
 & \Hom_{R^{\op}}(\Ext_R^n(V_0,R),Q).}$$
Since $\Ext_R^n(\lambda_{V_0}^U,R)$ is an isomorphism, we have that
$\Hom_{R^{\op}}(\Ext_R^n(\lambda_{V_0}^U,R),Q)$ is also an isomorphism, 
and hence so is $\Tor_n^R(Q,\lambda_{V_0}^U)$. Identifying $\Tor_n^R(Q,U)$ 
with its image in $\Tor_n^R(Q,V_0)$, that is, we have  
$$\Tor_n^R(Q,U)=\Tor_n^R(Q,V_0)\eqno{(3.7)}$$
for any non-zero submodule $U$ of $V_0$.

Set
$$\mathcal{U}:=\{\text{all non-zero submodules of}\ V_0\}\quad\text{and}\quad
S:=\bigcap_{U\in\mathcal{U}}U.$$
By \cite[Proposition 1 and Theorem~2]{Iw}, we have $\pd_{R^{\op}}Q=n$.
It follows that
\begin{align*}
\Tor_n^R(Q,S)&=\Tor_n^R(Q,\cap_{U\in\mathcal{U}}U)\\
& =\cap_{U\in\mathcal{U}}\Tor_n^R(Q,U)\quad\text{(by Proposition~\ref{prop-3.1})}\\
&=\Tor_n^R(Q,V_0)\quad\text{(by (3.7))}\\
&\cong\Hom_{R^{\op}}(\Ext_R^n(V_0,R),Q)\neq 0,\quad\text{(by Lemma~\ref{lem-2.3}(2))}
\end{align*}
and thus $S\neq 0$. Let $T$ be a non-zero submodule of $S$. Then $T\in\mathcal{U}$, and hence
$$S=\cap_{U\in\mathcal{U}}U\subseteq T\subseteq S.$$
Thus $T=S$, which implies that $S$ is simple. Consequently, we conclude that any non-zero
finitely generated submodule of $E^n(R)$ contains a simple submodule.

Now let $N$ be any non-zero submodule of $E^n(R)$. Choose
$0\neq x\in N$ and set $U:=Rx$. Then $U$ contains a simple submodule $S$. 
Since $0\neq S\subseteq U\subseteq N$ and $S\subseteq\Soc(E^n(R))$, we have
$$0\neq S\subseteq N\cap\Soc(E^n(R)).$$
It follows that $\Soc(E^n(R))$ is essential in $E^n(R)$. 
\end{proof}

As an application of Theorem \ref{thm-3.2}, we get its module version.

\begin{theorem}\label{thm-3.3}
Let $R$ be a Gorenstein ring with $\id_RR=\id_{R^{op}}R=n$. 
Then $E^n(M)$ has an essential socle for any $0\neq M\in\Mod R$ with $\id_RM=n$. 
\end{theorem}

\begin{proof}
Let $M$ be a non-zero left $R$-module with $\id_RM=n$. By \cite[Theorem~2.5]{M}, we can write
$$E^n(M)=\oplus_{i\in I}E_i,$$
where each $E_i$ is a non-zero indecomposable injective left $R$-module.
It follows from \cite[Theorem~1]{IM} that each $E_i$ is isomorphic to a direct summand of $E^n(R)$. 
Identifying $E_i$ with such a summand, we have
$$0\neq E_i\cap\Soc(E^n(R))=\Soc(E_i)$$
by \cite[p.242, Exercise 12(3)]{L}, and thus $E_i$ has a simple submodule $S_i$ by Theorem \ref{thm-3.2}. 
Since $E_i$ is indecomposable, we have $E_i=E^0(S_i)$, and so 
$$E^n(M)=\oplus_{i\in I}E^0(S_i).$$
Then it holds
\begin{align*}
\Soc E^n(M)&=\Soc\oplus_{i\in I}E^0(S_i)\\
&=\oplus_{i\in I}\Soc E^0(S_i)\quad\text{(by \cite[p.242, Exercise 12(6)]{L})}\\
&=\oplus_{i\in I}S_i.
\end{align*}
Since each $S_i$ is essential in $E^0(S_i)$, it follows from \cite[Example 3.38]{L} that
$\oplus_{i\in I}S_i$ is essential in $E^n(M)$.
\end{proof}

\section*{Acknowledgments}
The research was partially supported by NSFC (Grant Nos. 12371038, 12471036) and 
the Hubei Provincial Natural Science Foundation of China (Grant No.~2026AFA094).


\begin{thebibliography}{99}
\bibitem{B} H. Bass, {\it On the ubiquity of Gorenstein rings}, Math. Z. \textbf{82}(1) (1963), 8--28.

\bibitem{EJ} E. E. Enochs and O. M. G. Jenda, Relative Homological Algebra, Vol. 1,
the second revised and extended edition,
de Gruyter Expositions in Math. {\bf 30}, Walter de Gruyter GmbH \& Co. KG, Berlin, 2011.

\bibitem{FI}
K.R.~Fuller and Y.~Iwanaga,
\emph{On $n$-Gorenstein rings and Auslander rings of low injective dimension},
in \emph{Representations of Algebras} (Ottawa, ON, 1992),
CMS Conf. Proc. \textbf{14}, Amer. Math. Soc., Providence, RI,
1993, pp.175--183.


\bibitem{GN}
S.~Goto and K.~Nishida,
\emph{Towards a theory of Bass numbers with application to Gorenstein algebras},
Colloq. Math. \textbf{91}(2) (2002), 191--253.


\bibitem{Ho}
M.~Hoshino,
\emph{Noetherian rings of self-injective dimension two},
Comm. Algebra \textbf{21}(4) (1993), 1071--1094.


\bibitem{HC}
Z.Y.~Huang and F.C.~Cheng, 
\emph{On homological dimensions of simple modules over non-commutative rings}, Comm. Algebra 
\textbf{24}(10) (1996), 3259--3264.

\bibitem{HW}
Z.Y.~Huang and Y.~Wang,
\emph{The socle of the last term in a minimal injective resolution},
Sci. China Math. \textbf{53}(7) (2010), 1715--1721.

\bibitem{Iw}
Y.~Iwanaga,
\emph{On rings with finite self-injective dimension II},
Tsukuba J.~Math. \textbf{4}(1) (1980), 107--113.

\bibitem{IM}
Y.~Iwanaga and J.-I.~Miyachi,
\emph{Modules of the highest homological dimension over a Gorenstein ring},
in \emph{Trends in the Representation Theory of Finite-dimensional Algebras} (Seattle, WA, 1997), 
Contemp. Math. \textbf{229}, Amer. Math. Soc., Providence, RI, 1998, pp.193--199.

\bibitem{IS}
Y.~Iwanaga and H.~Sato,
\emph{On Auslander's $n$-Gorenstein rings},
J.~Pure Appl. Algebra \textbf{106}(1) (1996), 61--76.

\bibitem{L}
T.Y.~Lam,
\emph{Lectures on Modules and Rings},
Grad. Texts in Math. \textbf{189}, Springer-Verlag, New York, 1999.

\bibitem{M}
E. Matlis,
\emph{Injective modules over Noetherian rings}, Pacific J. Math. \textbf{8}(3) (1958), 511--528.


\bibitem{SZH}
W.L.~Song, X.J.~Zhang and Z.Y.~Huang,
\emph{The socle of the last term in the minimal injective resolution
of a Gorenstein module}, Osaka J.~Math. \textbf{56}(1) (2019), 123--132.

\end{thebibliography}
\end{document}